\documentclass[11pt,reqno]{amsart}
\usepackage{amsmath, amsthm, amscd, amssymb, amsfonts, latexsym, graphicx}
\usepackage{fullpage}
\usepackage{bm}
\usepackage[all]{xy}
\usepackage{tikz}
\usepackage{mathrsfs}
\usepackage{xcolor}
\usepackage[hidelinks]{hyperref}

\usetikzlibrary{patterns,decorations.pathreplacing}

\numberwithin{equation}{section}

\theoremstyle{plain}
\newtheorem{lem}{Lemma}[section]
\newtheorem{lemma}[lem]{Lemma}
\newtheorem{theorem}[lem]{Theorem}
\newtheorem{corollary}[lem]{Corollary}
\newtheorem{proposition}[lem]{Proposition}

\theoremstyle{definition}

\newtheorem{conjecture}[lem]{Conjecture}

\newcommand{\eps}{\varepsilon}
\newcommand{\cB}{\mathcal B}

\newcommand{\sym}{\operatorname{sym}}
\newcommand{\vol}{\operatorname{vol}}
\newcommand{\supp}{\operatorname{supp}}

\DeclareMathOperator{\Rea}{Re}
\newcommand{\HH}{\mathbb H}
\newcommand{\RR}{\mathbb R}
\newcommand{\CC}{\mathbb C}
\newcommand{\ZZ}{\mathbb Z}
\newcommand{\NN}{\mathbb N}
\newcommand{\ii}{\mathrm i}
\newcommand{\dd}{\,\mathrm d}

\newcommand{\PSL}{\mathrm{PSL}}

\newcommand{\Cc}{C_c^\infty}
\newcommand{\cW}{\mathcal W}

\newcommand{\mean}[1]{\frac{3}{\pi}\int_X #1\dd\mu}

\title[Conditional Equidistribution on the Imaginary Axis]
{Restricted Quantum unique ergodicity under Generalized Lindel\"of conjecture}
\author{Qingfeng Sun}
	\address{School of Mathematics and Statistics, Shandong University, Weihai\\Weihai, Shandong 264209, China}
	\address{State Key Laboratory of Cryptography and Digital Economy
Security \\ Shandong University, Weihai, 264209, China}
	\email{qfsun@sdu.edu.cn}
\author{Qizhi Zhang}
\address{School of Mathematics and Statistics, Shandong University, Weihai\\Weihai, Shandong 264209, China}
\email{qzzhang@mail.sdu.edu.cn}
\thanks{Q. Sun was partially supported by the National Natural Science Foundation of China
(Grant No. 12471005) and
the Natural Science Foundation of Shandong Province (Grant No. ZR2023MA003). Q. Zhang was supported by the Program of China Scholarship Council (Grant No. 202606220052).}
\date{}

\hypersetup{pdfauthor={Qizhi Zhang},pdftitle={Conditional Equidistribution on the Imaginary Axis via Shifted Sums and Mellin Separation}}

\begin{document}

\begin{abstract}
Assume the generalized Lindel\"of hypothesis. Let $u$ be an even
Hecke--Maass cusp form on $\PSL_2(\ZZ)\backslash\HH$, normalized in
$L^2$ with respect to hyperbolic area, with Laplace eigenvalue
$1/4+T^2$. For every fixed $\psi\in\Cc((0,\infty);\CC)$ and every
$\eps>0$, we prove that
\[
\int_0^\infty \psi(y)|u(\ii y)|^2\frac{\dd y}{y}
=\frac{6}{\pi}\int_0^\infty \psi(y)\frac{\dd y}{y}
+O_{\psi,\eps}(T^{-1/4+\eps})
\]
as $T\to\infty$, uniformly in $u$. In particular, the restriction mass
equidistributes on every fixed compact segment of the imaginary axis
in $\HH$.
\end{abstract}
\maketitle

\section{Introduction}
Let
\[
\Gamma=\PSL_2(\ZZ),\qquad
X=\Gamma\backslash\HH,\qquad
\dd\mu(z)=\frac{\dd x\dd y}{y^2},\qquad
\vol(X)=\frac{\pi}{3}.
\]
We use the nonnegative Laplacian
$\Delta=-y^2(\partial_x^2+\partial_y^2)$.
Throughout, $u$ is a real-valued even Hecke--Maass cusp form satisfying
\begin{equation}\label{eq:normalization}
\Delta u=\left(\frac14+T^2\right)u,\qquad
\int_X |u(z)|^2\dd\mu(z)=1,\qquad
u(-\overline z)=u(z),\qquad T\ge 2.
\end{equation}

The arithmetic quantum unique ergodicity (QUE) theorem of
Lindenstrauss and Soundararajan (\cite{Lindenstrauss} and \cite{Sound})
gives
\[
\int_X f(z)|u(z)|^2\dd\mu(z)
=\frac{3}{\pi}\int_X f(z)\dd\mu(z)+o_f(1)
\qquad (T\to\infty)
\]
for every fixed $f\in\Cc(X)$. Restricted QUE concerns the measures
obtained by restricting $|u|^2$ to a geodesic and integrating against
hyperbolic arclength. The geodesic has zero hyperbolic area, so the
surface asymptotic alone does not determine the restricted mass.
Young formulated the following conjecture for the geodesic joining
$0$ to $\ii\infty$ and proved its Eisenstein analogue
\cite[Conjecture~1.1 and Theorem~1.2]{Young16}.

\begin{conjecture}\label{conj:intro-rque}
As $u$ ranges over the forms in \eqref{eq:normalization}, for every
real-valued $\psi\in\Cc((0,\infty))$ one has
\begin{equation}\label{eq:intro-rque-prediction}
\lim_{T\to\infty}\int_0^\infty
\psi(y)|u(\ii y)|^2\frac{\dd y}{y}
=\frac{6}{\pi}\int_0^\infty\psi(y)\frac{\dd y}{y}.
\end{equation}
\end{conjecture}

Young later showed that the arithmetic position of a vertical
geodesic affects the Eisenstein restriction asymptotic and proposed
an analogous formula for Hecke--Maass cusp forms
\cite[Theorem~1.1 and Conjecture~1.3]{Young18}.
Unconditional density-one results on the imaginary axis are known
for holomorphic Hecke cusp forms by Zenz \cite[Corollary~1.4]{Zenz},
for holomorphic Hecke cusp forms in short weight intervals by Sun
and Zhang \cite[Corollary~1.4]{SZhol}, and for even Hecke--Maass cusp
forms in short spectral intervals by Sun and Zhang
\cite[Corollary~1.3]{SZ}.

In his 2024 lecture, Humphries \cite{Humphries} presented a period
method for effective restricted QUE on closed geodesics under the
generalized Lindel\"of hypothesis (GLH) and indicated that it also
applies to vertical geodesics from a rational point to $\ii\infty$.
For the imaginary axis we give a shifted-convolution proof with an
explicit power saving for each fixed smooth weight, uniform over
individual even forms.

\begin{theorem}\label{thm:main}
Assume the generalized Lindelof hypothesis \textup{(GLH)}.
For every fixed $\psi\in\Cc((0,\infty);\CC)$ and every $\eps>0$,
\begin{equation}\label{eq:main}
\int_0^\infty\psi(y)|u(iy)|^2\frac{\dd y}{y}
=\frac6\pi\int_0^\infty\psi(y)\frac{\dd y}{y}
+O_{\psi,\eps}(T^{-1/4+\eps})
\end{equation}
uniformly over the even Hecke--Maass cusp forms satisfying
\eqref{eq:normalization}.
\end{theorem}

The interval version of this result concerns the $L^2$ norm of a
restriction. For fixed $0<a<b<\infty$, let
$\gamma_{a,b}=\{iy:a\le y\le b\}\subset\HH$. Hyperbolic arclength
is $\dd s=\dd y/y$, and
\[
\|u|_{\gamma_{a,b}}\|_{L^2(\dd s)}^2
=\int_a^b|u(iy)|^2\frac{\dd y}{y}.
\]
The length of this segment is $\log(b/a)$. Bounds for these norms
have also played a role in the study of nodal domains of Maass forms.

Ghosh, Reznikov, and Sarnak \cite[Theorems~1.3 and~6.1]{GRS}
proved an upper bound of $O_\eps(T^\eps)$ for the squared
restriction norm on the imaginary axis, together with positive
lower bounds on sufficiently long fixed compact segments.
Jung \cite[Corollary~1.12]{Jung} obtained a positive lower bound
on every fixed compact segment for a density-one family of even
Hecke--Maass forms. The extension of this lower bound to all even
Hecke--Maass cusp forms follows from Jang and Jung
\cite[Corollary~3.2]{JJ}.
These bounds concern the size of the restriction norm, while the
restricted QUE asymptotic specifies its limiting value.

We record the corresponding interval formulation of
Theorem~\ref{thm:main}. It follows by approximating the
characteristic function of $[a,b]$ from above and below by smooth
functions and using the positivity of $|u|^2$.

\begin{corollary}\label{cor:interval}
Assume GLH. For every fixed $0<a<b<\infty$,
\begin{equation}\label{eq:interval}
\int_a^b|u(iy)|^2\frac{\dd y}{y}
=\frac6\pi\log(b/a)+o(1),
\end{equation}
as $T$ tends to infinity.
\end{corollary}

\paragraph{\textbf{Outline of proof.}}
The Fourier expansion reduces the restriction integral to a sum
over two Fourier indices. We separate the diagonal and split the
off-diagonal terms at shifts of size $R$. Poincar\'e series express
the diagonal and the main part of the short shifts as integrals
over $X$. Effective QUE under GLH supplies the diagonal main term
and bounds the short shifts.

The essential observation for the long shifts is that both Fourier
indices become small compared with $T$. The estimates of Sun and
Zhang \cite{SZ} allow us to restrict, with a negligible error, to
$n^2-m^2\ll T^{1+\beta}$, where $\beta>0$ can be chosen arbitrarily
small. Since $n-m>R$, we have
\[
m+n\ll\frac{T^{1+\beta}}R.
\]
For $R=T^\alpha$ and $\beta<\alpha$, both indices are therefore
$o(T)$. In this range the Bessel expansion in
Lemma~\ref{lem:low-kernel} produces the factors $(n/m)^{\pm\ii T}$.
After smooth dyadic decomposition, double Mellin inversion
separates the remaining weight and gives products
\[
L(s+\omega\ii T,u)L(z-\omega\ii T,u),\qquad \omega=\pm1.
\]

Moving both contours to the critical line, we apply GLH together
with the uniform Mellin estimate. Summing together gives the bound of long shifts.

\section{Preliminaries}

\subsection{Hecke--Maass forms}

Write $e(x)=e^{2\pi\ii x}$. For $u$ as in
\eqref{eq:normalization}, the Fourier expansion is
\begin{equation}\label{eq:Fourier}
u(x+iy)=\sqrt y\sum_{n\neq0}
\rho_u(n)K_{\ii T}(2\pi|n|y)e(nx).
\end{equation}
Evenness and the Hecke relations give
\[
\rho_u(n)=\rho_u(-n)=\rho_u(1)\lambda_u(n)
\quad(n\geq1),\qquad \lambda_u(1)=1.
\]
Put $L_u=L(1,\sym^2u)$. The $L^2$ normalization gives
\begin{equation}\label{eq:rho}
|\rho_u(1)|^2=\frac{2\cosh(\pi T)}{L_u};
\end{equation}
see \cite[p.~790]{LS04} or \cite[Section~2]{SZ}.
The bounds of Hoffstein--Lockhart \cite{HL} and Li \cite{Li},
together with \cite[Lemma~3.6]{GRS}, give
\begin{equation}\label{eq:coeff}
T^{-\eps}\ll_\eps L_u\ll_\eps T^\eps,\qquad
\sum_{n\leq Z}|\lambda_u(n)|^2
\ll_\eps Z(ZT)^\eps
\quad(Z\geq1).
\end{equation}
The Kim--Sarnak bound \cite{KS} is
\begin{equation*}
|\lambda_u(n)|\leq d(n)n^{7/64}\qquad(n\geq1),
\end{equation*}
where $d(n)$ is the divisor function.

The standard $L$-function attached to $u$ is
\begin{equation*}
L(s,u)=\sum_{n\geq1}\frac{\lambda_u(n)}{n^s}
=\prod_p\left(1-\lambda_u(p)p^{-s}+p^{-2s}\right)^{-1}
\qquad(\Rea s>1),
\end{equation*}
which can be analytically continued to the entire complex plane $\mathbb{C}$.
Under GLH, for every fixed $\nu>0$,
\begin{equation*}
\left|L\left(\tfrac12+\ii t,u\right)\right|
\ll_\nu(1+T+|t|)^\nu\qquad(t\in\RR),
\end{equation*}
uniformly in $u$ and $t$.

\subsection{Effective QUE under GLH}

Under GLH, we also have the following effective QUE estimate on $X$.
For $f\in\Cc(X)$, write
\[
J_u(f)=\int_X f(z)|u(z)|^2\dd\mu(z).
\]
With the nonnegative self-adjoint Laplacian $\Delta$ on $L^2(X)$, put
\begin{equation}\label{eq:Watson-norm}
\|f\|_{L^{2,s}(X)}=\|(1+\Delta)^{s/2}f\|_{L^2(X)},
\end{equation}
where fractional powers are defined by spectral calculus. For $f\in\Cc(X)$, spectral calculus
and integration by parts give
\begin{equation}\label{eq:Sobolev-energy}
\|f\|_{L^{2,1}(X)}^2
=\langle(1+\Delta)f,f\rangle_{L^2(X)}
=\|f\|_2^2+\int_X|\nabla_{\mathrm{hyp}}f|^2\dd\mu,
\end{equation}
where $|\nabla_{\mathrm{hyp}}f|^2
=y^2(|\partial_xf|^2+|\partial_yf|^2)$.
Taking $A=(1+\Delta)^{1/2}$, Cauchy--Schwarz gives
\begin{equation}\label{eq:Sobolev-interpolation}
\|f\|_{L^{2,1/2}(X)}^2
=\langle Af,f\rangle
\le\|Af\|_2\|f\|_2
=\|f\|_{L^{2,1}(X)}\|f\|_2.
\end{equation}

\begin{lemma}
\label{prop:Watson}
Assume GLH. For every fixed $\eta>0$ and every $f\in\Cc(X)$,
\begin{equation*}
\left|J_u(f)-\mean{f}\right|
\ll_{X,\eta} T^{-1/2+2\eta}\|f\|_{L^{2,1/2}(X)},
\end{equation*}
with an implied constant independent of $u$ and $f$.
\end{lemma}

\begin{proof}
Apply \cite[Section 4.2, Corollary 1]{Watson} under GLH with
$d_B=1$, $a=0$, and exponent $\eta>0$.
Watson's eigenvalue is $1/4+T^2$ and his Laplacian is $-\Delta$,
so his Sobolev norm agrees with \eqref{eq:Watson-norm} and
$(1/4+T^2)^{-1/4+\eta}\ll_\eta T^{-1/2+2\eta}$.
\end{proof}

\section{Bessel integrals}

Fix $\psi\in\Cc((0,\infty);\CC)$ and choose $0<a<b<\infty$
such that $\supp\psi\subset[a,b]$. For $m,n\ge1$, define
\begin{equation}\label{eq:Bdef}
\cB_{m,n}(T)
=16\cosh(\pi T)\int_0^\infty
\psi(y)K_{\ii T}(2\pi my)K_{\ii T}(2\pi ny)\dd y.
\end{equation}
We first record estimates that truncate sums involving $\cB_{m,n}$.
We then obtain an expansion of $\cB_{m,n}(T)$ when both Bessel
arguments are at most $T/2$. The factors $(n/m)^{\pm\ii T}$ in this
expansion will allow
us to use $L$-functions in Section~\ref{sec:Mellin}.

\subsection{Truncation}

The estimates of Sun and Zhang allow us to discard terms with
large $m$, $n$, or $|n^2-m^2|$, even after multiplying by a bounded weight.

\begin{lemma}
\label{lem:truncation}
There is a constant $C_\psi>0$ such that, for every fixed
$\beta,B>0$ and every bounded complex-valued weight
$b_{u,T}$ on $\NN^2$,
\begin{equation*}
\frac1{2L_u}
\sum_{(m,n)\notin\Omega_{\beta,T}}
\left|b_{u,T}(m,n)\lambda_u(m)\lambda_u(n)\cB_{m,n}\right|
\ll_{\psi,\beta,B}\|b_{u,T}\|_\infty T^{-B},
\end{equation*}
where
\[
\Omega_{\beta,T}
=\{(m,n)\in\NN^2:\ m,n\le C_\psi T,\quad
|n^2-m^2|\le T^{1+\beta}\}.
\]
\end{lemma}

\begin{proof}
The weighted statement follows by factoring out
$\|b_{u,T}\|_\infty$ and applying the tail estimates in
\cite[Lemmas~2.2 and~3.4]{SZ}, together with
\eqref{eq:coeff}.
\end{proof}

\subsection{Integral expression}

We begin with an asymptotic formula for $K_{\ii T}(x)$.

\begin{lemma}
\label{lem:WKB}
For $T\ge1$ and $0<x\le3T/4$, put
\[
k_T(x)=\sqrt{T^2-x^2},
\qquad
\Theta_T(x)=T\operatorname{arcosh}(T/x)-k_T(x)-\frac\pi4.
\]
Then
\begin{equation*}
K_{\ii T}(x)=\sqrt{2\pi}\,e^{-\pi T/2}k_T(x)^{-1/2}
\left\{\cos\Theta_T(x)+O(T^{-1})\right\},
\end{equation*}
with an absolute implied constant.
\end{lemma}

\begin{proof}
Put $q=x/T$ and $\Theta=\Theta_T(x)$.
By Andersen--Duke \cite[Proposition A.1]{AD22}, with their
$\nu=T$ and $z=q$, uniformly for $T\ge1$ and $0<q\le3/4$,
\[
e^{\pi T/2}K_{\ii T}(Tq)
=\frac{\sqrt{2\pi}}{T^{1/2}(1-q^2)^{1/4}}\left\{
\cos\Theta+
\frac{3q^2+2}{24T(1-q^2)^{3/2}}\sin\Theta\right\}
+O(T^{-5/2}).
\]
This gives the claimed asymptotic formula.
\end{proof}

Define
\begin{equation*}
\cW_\psi(v)=\int_0^\infty
\psi(y)e^{-\ii\pi^2vy^2}\dd y.
\end{equation*}
With $s=y^2$, the function
$g_\psi(s)=\psi(\sqrt s)/(2\sqrt s)$ is also smooth on $(0,\infty)$,
which is supported in $[a^2,b^2]$. Its zero extension belongs to
$\Cc(\RR)$, and
\[
\cW_\psi(v)=\int_\RR g_\psi(s)e^{-\ii\pi^2vs}\dd s.
\]
It's not hard to show that $\cW_\psi$ is a Schwartz function on $\RR$.

Applying Lemma~\ref{lem:WKB} inside the integral defining $\cB_{m,n}(T)$
and expanding the phase difference isolates the factors
$(n/m)^{\pm\ii T}$. The following lemma gives the resulting formula.

\begin{lemma}\label{lem:low-kernel}
Suppose $m<n$ and $2\pi bn\leq T/2$.
Put $\Delta_{m,n}=n^2-m^2$.
Then
\begin{align*}
\cB_{m,n}(T)
={}&\frac{4\pi}{T}
\sum_{\omega\in\{-1,1\}}
(n/m)^{\omega\ii T}
\cW_\psi(\omega\Delta_{m,n}/T)
\\
&+O_\psi\left(
T^{-2}+
\frac{m^2+n^2}{T^3}
\left(1+\frac{\Delta_{m,n}}{T}\right)
\right).
\end{align*}
For arbitrary positive integers $m,n$ with
$2\pi b\max(m,n)\le T/2$, including $m=n$, one also has
\begin{equation*}
\cB_{m,n}(T)\ll_\psi T^{-1}.
\end{equation*}
\end{lemma}

\begin{proof}
For $y\in[a,b]$, write
\[
x_m(y)=2\pi my,\qquad x_n(y)=2\pi ny,\qquad
k_m(y)=\sqrt{T^2-x_m(y)^2},\qquad
k_n(y)=\sqrt{T^2-x_n(y)^2}.
\]
The assumption $2\pi bn\le T/2$ gives $k_m(y),k_n(y)\ge\sqrt3T/2$.
Substituting the asymptotic formula in Lemma~\ref{lem:WKB} into the definition of
$\cB_{m,n}(T)$ therefore gives
\begin{align*}
\cB_{m,n}(T)
={}&16\pi(1+e^{-2\pi T})
\int_a^b\frac{\psi(y)}{\sqrt{k_m(y)k_n(y)}}\\
&\hspace{3em}\times
\cos\Theta_T(x_m(y))
\cos\Theta_T(x_n(y))\dd y
+O_\psi(T^{-2}).
\end{align*}
Indeed, the error in each cosine factor is $O(T^{-1})$,
whereas $(k_mk_n)^{-1/2}\ll T^{-1}$.
Set
\[
D(y)=\Theta_T(x_m(y))-\Theta_T(x_n(y)),\qquad
S(y)=\Theta_T(x_m(y))+\Theta_T(x_n(y)).
\]
Using $2\cos u\cos v=\cos(u-v)+\cos(u+v)$, we obtain
\begin{align*}
\cB_{m,n}(T)
={}&8\pi\int_a^b
\frac{\psi(y)\cos D(y)}{\sqrt{k_m(y)k_n(y)}}\dd y\\
&+8\pi\int_a^b
\frac{\psi(y)\cos S(y)}{\sqrt{k_m(y)k_n(y)}}\dd y
+O_\psi(T^{-2}).
\end{align*}

We first estimate the integral containing $\cos S(y)$. Since
$\Theta_T'(x)=-k_T(x)/x$,
\[
-S'(y)=\frac{k_m(y)+k_n(y)}{y}\asymp_\psi T,
\qquad
k_j'(y)=-\frac{(2\pi j)^2y}{k_j(y)}=O_\psi(T)
\quad(j=m,n).
\]
Consequently,
\[
\frac{\dd}{\dd y}
\left(\frac{y\psi(y)}{\sqrt{k_m(y)k_n(y)}\,(k_m(y)+k_n(y))}\right)
=O_\psi(T^{-2}).
\]
Integration by parts now yields
\begin{align*}
\int_a^b\frac{\psi(y)\cos S(y)}{\sqrt{k_m(y)k_n(y)}}\dd y
&=\int_a^b\sin S(y)\,
\frac{\dd}{\dd y}
\left(\frac{y\psi(y)}{\sqrt{k_m(y)k_n(y)}\,(k_m(y)+k_n(y))}\right)\dd y\\
&=O_\psi(T^{-2}).
\end{align*}
Thus
\[
\cB_{m,n}(T)
=8\pi\int_a^b\frac{\psi(y)\cos D(y)}{\sqrt{k_m(y)k_n(y)}}\dd y
+O_\psi(T^{-2}).
\]

For $D(y)$, the identity
\[
D(y)=\int_{x_m(y)}^{x_n(y)}\frac{k_T(v)}v\dd v
\]
and the uniform expansion
$k_T(v)=T-v^2/(2T)+O(v^4/T^3)$ for $0<v\le T/2$ give
\begin{align*}
D(y)
&=T\log\frac{x_n(y)}{x_m(y)}
-\frac{x_n(y)^2-x_m(y)^2}{4T}
+O\left(\frac{x_n(y)^4-x_m(y)^4}{T^3}\right)\\
&=T\log(n/m)-\frac{\pi^2\Delta_{m,n}y^2}{T}
+O_\psi\left(\frac{\Delta_{m,n}(m^2+n^2)}{T^3}\right).
\end{align*}
Also,
\begin{align*}
\frac1{\sqrt{k_m(y)k_n(y)}}
&=\frac1T
\left(1-\frac{x_m(y)^2}{T^2}\right)^{-1/4}
\left(1-\frac{x_n(y)^2}{T^2}\right)^{-1/4}\\
&=\frac1T+O_\psi\left(\frac{m^2+n^2}{T^3}\right).
\end{align*}
By the mean value theorem,
\[
\cos D(y)
=\cos\left(T\log(n/m)-\frac{\pi^2\Delta_{m,n}y^2}{T}\right)
+O_\psi\left(\frac{\Delta_{m,n}(m^2+n^2)}{T^3}\right).
\]
Combining these estimates gives
\begin{align*}
\cB_{m,n}(T)
={}&\frac{8\pi}{T}\int_a^b\psi(y)
\cos\left(T\log(n/m)-\frac{\pi^2\Delta_{m,n}y^2}{T}\right)\dd y\\
&+O_\psi\left(T^{-2}+\frac{m^2+n^2}{T^3}
+\frac{\Delta_{m,n}(m^2+n^2)}{T^4}\right).
\end{align*}
Finally, the main term is exactly
\begin{align*}
&\frac{4\pi}{T}\sum_{\omega\in\{-1,1\}}
(n/m)^{\omega\ii T}
\int_a^b\psi(y)e^{-\omega\ii\pi^2\Delta_{m,n}y^2/T}\dd y\\
&\hspace{2em}=
\frac{4\pi}{T}\sum_{\omega\in\{-1,1\}}
(n/m)^{\omega\ii T}\cW_\psi(\omega\Delta_{m,n}/T),
\end{align*}
which proves the claimed expansion.

For the final bound, apply Lemma~\ref{lem:WKB} separately to both
indices to obtain
\[
e^{\pi T/2}|K_{\ii T}(2\pi jy)|\ll T^{-1/2}
\qquad(j=m,n,\ y\in[a,b]).
\]
Substitution into the defining integral for $\cB_{m,n}(T)$
gives $\cB_{m,n}(T)\ll_\psi T^{-1}$, also when $m=n$.
\end{proof}

\section{Decomposition of the restriction integral}

For the fixed test function $\psi$, write
\begin{equation}\label{eq:restriction}
I_u(\psi)=\int_0^\infty\psi(y)|u(iy)|^2\frac{\dd y}{y},
\qquad
M(\psi)=\frac6\pi\int_0^\infty\psi(y)\frac{\dd y}{y}.
\end{equation}
Substituting the Fourier expansion \eqref{eq:Fourier} into \eqref{eq:restriction} and
using \eqref{eq:rho} and \eqref{eq:Bdef}, we obtain
\begin{align}
I_u(\psi)
&=4|\rho_u(1)|^2\sum_{m,n\ge1}\lambda_u(m)\lambda_u(n)
\int_0^\infty\psi(y)K_{\ii T}(2\pi my)K_{\ii T}(2\pi ny)\dd y
\nonumber\\
&=\frac1{2L_u}\sum_{m,n\ge1}\lambda_u(m)\lambda_u(n)\cB_{m,n}(T).
\label{eq:I-sum}
\end{align}

We first split the sum in \eqref{eq:I-sum} into the diagonal and
positive shifts. Put
\begin{equation}\label{eq:diagonal-shifts}
D_u=\frac{1}{2L_u}\sum_{n\geq1}
\lambda_u(n)^2\cB_{n,n},\qquad
S_u(r)=\frac{1}{L_u}\sum_{n\geq1}
\lambda_u(n)\lambda_u(n+r)\cB_{n,n+r}.
\end{equation}
Choose once and for all a smooth function $\chi$ on $\RR$
with values in $[0,1]$ such that
\begin{equation*}
\chi(t)=1\quad(t\leq1),\qquad
\chi(t)=0\quad(t\geq2).
\end{equation*}
For $R\geq1$, we then split the off-diagonal sum into
\begin{align}
\mathcal O_{\mathrm{short}}
&=\sum_{r\geq1}\chi(r/R)S_u(r),\nonumber\\
\mathcal O_{\mathrm{long}}
&=\sum_{r\geq1}\bigl(1-\chi(r/R)\bigr)S_u(r).
\label{eq:shift-split}
\end{align}
Thus
\begin{equation}\label{eq:mass-decomposition}
I_u(\psi)=D_u+\sum_{r\geq1}S_u(r)
=D_u+\mathcal O_{\mathrm{short}}+\mathcal O_{\mathrm{long}}.
\end{equation}

\subsection{Poincar\'e series}
To relate $D_u$ and $S_u(r)$ to integrals over $X$, let
$\Gamma_\infty$ be the stabilizer of the cusp at infinity.
For $r\in\ZZ$, the Poincar\'e series $P_r(z)$ is defined by
\begin{equation}\label{eq:Pdef}
P_r(z)=
\sum_{\gamma\in\Gamma_\infty\backslash\Gamma}
\Im(\gamma z)\psi(\Im(\gamma z))e(r\Re(\gamma z)).
\end{equation}
To apply Lemma \ref{prop:Watson}, we therefore estimate the Sobolev norm of arbitrary finite linear combinations of Poincar\'e series.

\begin{lemma}\label{lem:Poincare}
The series in \eqref{eq:Pdef} are locally finite. Their sums belong
to $\Cc(X)$, and their supports lie in a common compact set
$K_\psi$, independent of $r$.
For every finitely supported sequence $(c_r)_{r\in\ZZ}$,
\begin{equation*}
\left\|\sum_r c_rP_r\right\|_{L^{2,1/2}(X)}
\ll_\psi
\left(\sum_r|c_r|^2\right)^{1/4}
\left(\sum_r(1+r^2)|c_r|^2\right)^{1/4}.
\end{equation*}
\end{lemma}

\begin{proof}
Choose a nonnegative $h\in\Cc((0,\infty))$ with
$h(y)\ge y|\psi(y)|$ and $h>0$ on $[a,b]$.
By \cite[Section~1.8]{Garrett}, the incomplete Eisenstein series
\[
E_h(z)=\sum_{\gamma\in\Gamma_\infty\backslash\Gamma}
 h(\Im(\gamma z))
\]
is locally finite and belongs to $\Cc(X)$.
It follows that \eqref{eq:Pdef} is locally finite, and
$|P_r|\le E_h$ gives $\supp P_r\subset K_\psi:=\supp E_h$.
Moreover, writing $c_h=\min_{[a,b]}h>0$, we have the uniform bound
\[
N_{a,b}(z):=\#\{\gamma\in\Gamma_\infty\backslash\Gamma:
a\le\Im(\gamma z)\le b\}
\le c_h^{-1}E_h(z)\ll_\psi1.
\]

Put
\[
F(x)=\sum_r c_re(rx),\qquad
g(x+iy)=y\psi(y)F(x),\qquad
P=\sum_r c_rP_r=\sum_\gamma g(\gamma z).
\]
Cauchy--Schwarz followed by integration over $X$ and unfolding to
$\Gamma_\infty\backslash\HH=\{0\le x<1,\ y>0\}$ gives
\begin{align*}
\|P\|_2^2
&\ll_\psi\int_0^\infty\int_0^1|g(x+iy)|^2
\frac{\dd x\dd y}{y^2}
\ll_\psi\int_0^1|F(x)|^2\dd x,\\
\|\nabla_{\mathrm{hyp}}P\|_2^2
&\ll_\psi\int_0^\infty\int_0^1
\bigl(|\partial_xg|^2+|\partial_yg|^2\bigr)\dd x\dd y\\
&\ll_\psi\int_0^1\bigl(|F(x)|^2+|F'(x)|^2\bigr)\dd x,
\end{align*}
where $\partial_xg=y\psi(y)F'(x)$ and
$\partial_yg=(\psi(y)+y\psi'(y))F(x)$.
Fourier orthogonality gives
\[
\int_0^1|F(x)|^2\dd x=\sum_r|c_r|^2,\qquad
\int_0^1|F'(x)|^2\dd x=4\pi^2\sum_r r^2|c_r|^2.
\]
Consequently, \eqref{eq:Sobolev-energy} implies
\[
\|P\|_2^2\ll_\psi\sum_r|c_r|^2,\qquad
\|P\|_{L^{2,1}(X)}^2\ll_\psi\sum_r(1+r^2)|c_r|^2.
\]
Applying \eqref{eq:Sobolev-interpolation} proves the claimed bound.
\end{proof}

Unfolding against the constant function gives
\begin{equation}\label{eq:Pmean}
\int_XP_r\dd\mu
=
\begin{cases}
\displaystyle\int_0^\infty\psi(y)\frac{\dd y}{y},&r=0,\\
0,&r\neq0.
\end{cases}
\end{equation}
The next lemma expresses $D_u$ and $S_u(r)$ in terms of $J_u(P_r)$.

\begin{lemma}\label{lem:periods}
One has
\begin{equation*}
D_u=2J_u(P_0)
\end{equation*}
and, for every integer $r\geq1$,
\begin{equation*}
S_u(r)=4J_u(P_r)
-\frac{1}{2L_u}\sum_{a'=1}^{r-1}
\lambda_u(a')\lambda_u(r-a')\cB_{a',r-a'}.
\end{equation*}
\end{lemma}

\begin{proof}
Unfolding gives
\[
J_u(P_r)=
\int_0^\infty\int_0^1
e(rx)\psi(y)|u(x+iy)|^2\dd x\,\frac{\dd y}{y}.
\]
For $r=0$, the $x$ integral pairs equal nonzero Fourier indices.
Their positive and negative parts contribute equally.
Substitution of \eqref{eq:rho} gives
\[
J_u(P_0)=\frac{1}{4L_u}
\sum_{n\geq1}\lambda_u(n)^2\cB_{n,n},
\]
which proves the identity for $D_u$.

For $r>0$, the Fourier indices are $k$ and $k+r$, with
$k\neq0,-r$.
The ranges $k\geq1$ and $k\leq-r-1$ give two identical infinite
sums. The remaining indices are $k=-a'$,
$1\leq a'\leq r-1$. Therefore
\[
J_u(P_r)=\frac{1}{4L_u}
\sum_{n\geq1}\lambda_u(n)\lambda_u(n+r)\cB_{n,n+r}
+\frac{1}{8L_u}
\sum_{a'=1}^{r-1}
\lambda_u(a')\lambda_u(r-a')\cB_{a',r-a'}.
\]
Rearranging gives the identity for $S_u(r)$.
\end{proof}

Put $\mathcal P_R=\sum_{r\ge1}\chi(r/R)P_r$.
Multiplying the identity for $S_u(r)$ in Lemma~\ref{lem:periods} by $\chi(r/R)$ and summing
over $r\ge1$ gives
\begin{equation}\label{eq:short-decomposition}
\mathcal O_{\mathrm{short}}
=\mathcal O_{\mathrm{short}}^{(1)}
+\mathcal O_{\mathrm{short}}^{(2)},
\end{equation}
where
\begin{align}
\mathcal O_{\mathrm{short}}^{(1)}
&=4\sum_{r\ge1}\chi(r/R)J_u(P_r)
=4J_u(\mathcal P_R),\nonumber\\
\mathcal O_{\mathrm{short}}^{(2)}
&=-\frac1{2L_u}\sum_{r\ge1}\chi(r/R)
\sum_{a'=1}^{r-1}\lambda_u(a')\lambda_u(r-a')\cB_{a',r-a'}.
\label{eq:short-components}
\end{align}
The second term subtracts the finite sum corresponding to
$-r<k<0$ in the preceding calculation.

\subsection{Short shifts}

We will apply effective QUE to $J_u(\mathcal P_R)$ and estimate this
finite sum directly.
We first estimate $D_u$ using effective QUE.
\begin{lemma}\label{lem:diagonal-bound}
Assume GLH. For every fixed $\eta>0$,
\[
D_u=M(\psi)+O_{\psi,\eta}(T^{-1/2+2\eta}).
\]
\end{lemma}

\begin{proof}
Taking $c_0=1$ and $c_r=0$ for $r\ne0$ in
Lemma~\ref{lem:Poincare} gives $\|P_0\|_{L^{2,1/2}(X)}\ll_\psi1$.
Together with the identity for $D_u$ in Lemma~\ref{lem:periods},
\eqref{eq:Pmean}, and Lemma~\ref{prop:Watson}, this proves the result.
\end{proof}

We next bound $\mathcal O_{\mathrm{short}}^{(1)}$ by applying
effective QUE to the sum $\mathcal P_R$ of the test functions.

\begin{lemma}\label{lem:short-period-bound}
Assume GLH. For every $R\ge1$ and every fixed $\eta>0$,
\[
\left|\mathcal O_{\mathrm{short}}^{(1)}\right|
\ll_{\psi,\eta}T^{-1/2+2\eta}R.
\]
\end{lemma}

\begin{proof}
By \eqref{eq:short-components},
$\mathcal O_{\mathrm{short}}^{(1)}=4J_u(\mathcal P_R)$.
The mean of $\mathcal P_R$ is zero by \eqref{eq:Pmean}.
Since $0\le\chi\le1$ and $\chi(r/R)=0$ for $r\ge2R$,
Lemma~\ref{lem:Poincare} gives
\[
\|\mathcal P_R\|_{L^{2,1/2}(X)}
\ll_\psi
\left(\sum_{1\le r<2R}\chi(r/R)^2\right)^{1/4}
\left(\sum_{1\le r<2R}(1+r^2)\chi(r/R)^2\right)^{1/4}
\ll_\psi R.
\]
Applying Lemma~\ref{prop:Watson} to $\mathcal P_R$ proves the lemma.
\end{proof}

For $\mathcal O_{\mathrm{short}}^{(2)}$, we use the final bound in Lemma~\ref{lem:low-kernel}
and Cauchy--Schwarz to estimate the finite sum directly.

\begin{lemma}\label{lem:short-correction-bound}
If $R=T^\alpha$ with fixed $0<\alpha<1$, then for every fixed
$\eps>0$,
\[
\left|\mathcal O_{\mathrm{short}}^{(2)}\right|
\ll_{\psi,\alpha,\eps}T^{-1+\eps}R^2.
\]
\end{lemma}

\begin{proof}
Lemma~\ref{lem:low-kernel} implies
$\cB_{a',r-a'}\ll_\psi T^{-1}$ for $r\leq2R$,
once $T$ is sufficiently large in terms of $\psi$ and $\alpha$:
both indices are at most $2R=o(T)$.
By \eqref{eq:short-components} and the bound from Lemma~\ref{lem:low-kernel},
\begin{align*}
\left|\mathcal O_{\mathrm{short}}^{(2)}\right|
&\le\frac1{2L_u}\sum_{r\geq1}\chi(r/R)
\sum_{a'=1}^{r-1}
\left|\lambda_u(a')\lambda_u(r-a')\cB_{a',r-a'}\right|\\
&\ll_\psi T^{-1}L_u^{-1}
\left(\sum_{n\leq2R}|\lambda_u(n)|\right)^2.
\end{align*}
In the second line, we used $0\le\chi\le1$ and enlarged the sum
to all pairs with both indices at most $2R$.
Now put $\delta=\eps/4$. Cauchy--Schwarz gives the first inequality in
\[
\left(\sum_{n\leq2R}|\lambda_u(n)|\right)^2
\le2R\sum_{n\leq2R}|\lambda_u(n)|^2
\ll_\delta R^2(RT)^\delta,
\]
and the last bound follows from \eqref{eq:coeff} with exponent $\delta$.
Using also $L_u^{-1}\ll_\delta T^\delta$ from \eqref{eq:coeff}, we obtain
\[
\left|\mathcal O_{\mathrm{short}}^{(2)}\right|
\ll_{\psi,\delta}T^{-1+\delta}R^2(RT)^\delta
=T^{-1+(2+\alpha)\delta}R^2
\ll_{\psi,\alpha,\eps}T^{-1+\eps}R^2,
\]
since $(2+\alpha)\delta<\eps$ for $\alpha<1$.
\end{proof}

\section{Reduction of the long-shift sum}

Our aim in this section is to express the long-shift sum in a form
that can be estimated by Mellin inversion and GLH in
Section~\ref{sec:Mellin}. We first restrict $n^2-m^2$ and use the
lower bound on $n-m$ to bring both indices into the range of
Lemma~\ref{lem:low-kernel}. Applying that lemma gives the sums
$\mathcal H_\omega$ in Proposition~\ref{prop:long-reduction}; we also
bound the errors from the truncation and the
approximation in Lemma~\ref{lem:low-kernel}.

By \eqref{eq:diagonal-shifts} and \eqref{eq:shift-split}, we have
\[
\mathcal O_{\mathrm{long}}
=\frac1{L_u}\sum_{m<n}
\lambda_u(m)\lambda_u(n)
\bigl(1-\chi((n-m)/R)\bigr)\cB_{m,n}.
\]
To restrict $n^2-m^2$, for $R\geq1$ and $M>0$ define
\[
A_{m,n}=\bigl(1-\chi((n-m)/R)\bigr)
\chi((n^2-m^2)/M)\qquad(m<n),
\]
and set $A_{m,n}=0$ for $n\le m$.

\begin{lemma}\label{lem:long-local}
Let $\alpha$ and $\beta$ be fixed with
\begin{equation*}
0<\alpha\le\frac14,\qquad
0<\beta<\frac{\alpha}{4},
\end{equation*}
and set
\begin{equation*}
R=T^\alpha,\qquad M=T^{1+\beta},\qquad Y=\frac{2M}{R}.
\end{equation*}
Then, for every fixed $B>0$,
\begin{equation*}
\mathcal O_{\mathrm{long}}
=\frac1{L_u}\sum_{m<n}
\lambda_u(m)\lambda_u(n)
A_{m,n}\cB_{m,n}
+O_{\psi,\alpha,\beta,B}(T^{-B}).
\end{equation*}
If $A_{m,n}\neq0$, then
\begin{equation*}
n-m>R,\qquad n^2-m^2<2M,\qquad m+n<Y.
\end{equation*}
In particular, $m,n=o(T)$.
\end{lemma}

\begin{proof}
Apply Lemma~\ref{lem:truncation} with the present $\beta$ and
\[
b_{u,T}(m,n)=\mathbf1_{m<n}
\bigl(1-\chi((n-m)/R)\bigr)
\bigl(1-\chi((n^2-m^2)/M)\bigr).
\]
This weight has sup norm at most one and vanishes when
$|n^2-m^2|\le M=T^{1+\beta}$. Hence its contribution is
$O_{\psi,\beta,B}(T^{-B})$, proving the asserted identity.

The first two support conditions in the statement follow
from the cutoffs. Their quotient gives
\[
m+n=\frac{n^2-m^2}{n-m}<\frac{2M}{R}=Y.
\]
Since $Y/T=2T^{\beta-\alpha}=o(1)$, the last assertion follows.
\end{proof}

Keep $\alpha,\beta,R,M,Y$ as in Lemma~\ref{lem:long-local}.
The following bound controls the sum of the errors in the Bessel
expansion.

\begin{lemma}
\label{lem:band-mass}
Suppose $1\leq Y_0\leq M_0\leq T^2$.
Then
\begin{equation*}
\sum_{\substack{1\leq m<n\leq Y_0\\
n^2-m^2\leq2M_0}}
|\lambda_u(m)\lambda_u(n)|
\ll_\eps M_0T^\eps.
\end{equation*}
\end{lemma}

\begin{proof}
Fix $\eps>0$ and put $\delta=\eps/4$.
For $n\leq2\sqrt{M_0}$, Cauchy--Schwarz and
\eqref{eq:coeff} with exponent $\delta$ give
\begin{align*}
\sum_{m,n\leq2\sqrt{M_0}}
|\lambda_u(m)\lambda_u(n)|
&=\left(\sum_{n\leq2\sqrt{M_0}}|\lambda_u(n)|\right)^2\\
&\leq2\sqrt{M_0}\sum_{n\leq2\sqrt{M_0}}|\lambda_u(n)|^2\\
&\ll_\delta M_0(T\sqrt{M_0})^\delta
\leq M_0T^{2\delta},
\end{align*}
since $M_0\leq T^2$.

For the remaining indices, use dyadic intervals $Z\leq n<2Z$.
Let $\mathcal E_Z$ denote the pairs in the sum with
$Z\leq n<2Z$ and $n>2\sqrt{M_0}$.
The inequality $n^2-m^2\le2M_0$ gives $m^2\geq n^2-2M_0>n^2/2$; hence
\[
\frac{Z}{\sqrt2}<m<n<2Z,
\qquad
1\leq n-m=\frac{n^2-m^2}{n+m}\leq\frac{2M_0}{Z}.
\]
For each fixed $n$, there are $O(1+M_0/Z)$ possible values of $m$;
the same bound holds for $n$ when $m$ is fixed. If $\mathcal E_Z$
is nonempty, then $Z\leq Y_0\leq M_0$, so each of these counts is
$O(M_0/Z)$.
It follows that
\begin{align*}
\sum_{(m,n)\in\mathcal E_Z}|\lambda_u(m)\lambda_u(n)|
&\leq\frac12\sum_{(m,n)\in\mathcal E_Z}
\bigl(|\lambda_u(m)|^2+|\lambda_u(n)|^2\bigr)\\
&\ll\frac{M_0}{Z}\sum_{k\leq2Z}|\lambda_u(k)|^2\\
&\ll_\delta\frac{M_0}{Z}\,Z(ZT)^\delta
\leq M_0T^{3\delta},
\end{align*}
where the last step uses $Z\leq M_0\leq T^2$.
There are $O(\log(2T))$ dyadic ranges. Combining their contributions
with the first range and using $\log(2T)\ll_\delta T^\delta$ gives
\[
\sum_{\substack{1\leq m<n\leq Y_0\\n^2-m^2\leq2M_0}}
|\lambda_u(m)\lambda_u(n)|
\ll_\delta M_0T^{2\delta}
+M_0T^{3\delta}\log(2T)
\ll_\eps M_0T^{4\delta}=M_0T^\eps.
\]
\end{proof}

Now the range in Lemma~\ref{lem:long-local}
is already sufficient for applying Lemma~\ref{lem:low-kernel}.
Indeed, whenever $A_{m,n}\ne0$, we have
$n<Y=2T^{1+\beta-\alpha}=o(T)$, so $2\pi bn\le T/2$ for
sufficiently large $T$. Inserting the expansion from Lemma~\ref{lem:low-kernel} into
the identity in Lemma~\ref{lem:long-local} therefore produces the factors
$(n/m)^{\pm\ii T}$ in the sum. Their form provides the link to
$L$-functions: since
$(n/m)^{\omega\ii T}=m^{-\omega\ii T}n^{\omega\ii T}$,
Mellin inversion of the remaining weight separates the two
summation variables and leads to two $L$-functions with shifts
$\pm\ii T$. We carry out this step in Section~\ref{sec:Mellin}.
The following proposition records the sums obtained by substituting
the expansion from Lemma~\ref{lem:low-kernel} into the identity in
Lemma~\ref{lem:long-local} and bounds the resulting error.

\begin{proposition}
\label{prop:long-reduction}
For $\omega\in\{-1,1\}$, define
\begin{equation}\label{eq:Hsign}
\mathcal H_\omega
=\frac1{TL_u}\sum_{m,n\geq1}
\lambda_u(m)\lambda_u(n)(n/m)^{\omega\ii T}
A_{m,n}\cW_\psi\left(\omega\frac{n^2-m^2}{T}\right).
\end{equation}
Then, for every $\eps>0$,
\begin{equation*}
\mathcal O_{\mathrm{long}}
=4\pi\sum_{\omega=\pm1}\mathcal H_\omega
+O_{\psi,\alpha,\beta,\eps}
\left(
T^{-1+\beta+\eps}
+T^{-2\alpha+4\beta+\eps}
+T^{-2}
\right).
\end{equation*}
\end{proposition}

\begin{proof}
Take $B=2$ in Lemma~\ref{lem:long-local}. Whenever $A_{m,n}\ne0$,
the support conditions in that lemma give $n<Y=o(T)$, so
Lemma~\ref{lem:low-kernel} applies for sufficiently large $T$.
Substituting the expansion from Lemma~\ref{lem:low-kernel} into
the identity in Lemma~\ref{lem:long-local} gives
\begin{align*}
\mathcal O_{\mathrm{long}}
={}&\frac{4\pi}{TL_u}\sum_{\omega=\pm1}\sum_{m<n}
\lambda_u(m)\lambda_u(n)(n/m)^{\omega\ii T}
A_{m,n}\cW_\psi\left(\omega\frac{n^2-m^2}{T}\right)\\
&+O_\psi(E_1+E_2+E_3)
+O_{\psi,\alpha,\beta}(T^{-2}),
\end{align*}
where
\begin{align*}
E_1&=\frac1{T^2L_u}\sum_{m<n}
|\lambda_u(m)\lambda_u(n)|A_{m,n},\\
E_2&=\frac1{T^3L_u}\sum_{m<n}
|\lambda_u(m)\lambda_u(n)|A_{m,n}(m^2+n^2),\\
E_3&=\frac1{T^4L_u}\sum_{m<n}
|\lambda_u(m)\lambda_u(n)|A_{m,n}(m^2+n^2)(n^2-m^2).
\end{align*}
The quantities $E_1,E_2,E_3$ bound the contributions of the three
terms in the error in Lemma~\ref{lem:low-kernel}; the last
$O(T^{-2})$ term comes from Lemma~\ref{lem:long-local}.
Since $A_{m,n}=0$ for $n\le m$, the first term in the decomposition
is exactly $4\pi\sum_{\omega=\pm1}\mathcal H_\omega$ by
\eqref{eq:Hsign}.

To estimate $E_1,E_2,E_3$, we first bound the sum of
$|\lambda_u(m)\lambda_u(n)|A_{m,n}$. For sufficiently large $T$,
the parameter conditions in Lemma~\ref{lem:long-local} give
\[
1\leq Y\leq M=T^{1+\beta}\leq T^2.
\]
Thus Lemma~\ref{lem:band-mass} applies with $Y_0=Y$ and $M_0=M$.
Using $0\leq A_{m,n}\leq1$ and the support conditions
in Lemma~\ref{lem:long-local}, we apply Lemma~\ref{lem:band-mass} and
\eqref{eq:coeff}, each with exponent $\eps/2$, to obtain
\begin{align*}
\frac1{L_u}\sum_{m<n}|\lambda_u(m)\lambda_u(n)|A_{m,n}
\ll_\eps T^{\eps/2}
\sum_{\substack{1\leq m<n\leq Y\\n^2-m^2\leq2M}}
|\lambda_u(m)\lambda_u(n)|\ll_\eps  MT^\eps.
\end{align*}
It follows directly that
\[
E_1\ll_\eps\frac{MT^\eps}{T^2}
=T^{-1+\beta+\eps}.
\]

For $E_2$, we use $m^2+n^2\leq(m+n)^2<Y^2$ whenever
$A_{m,n}\ne0$. The same coefficient estimate therefore gives
\begin{align*}
E_2&\ll_\eps\frac{Y^2}{T^3}\,MT^\eps\ll T^{-2\alpha+3\beta+\eps}.
\end{align*}
For $E_3$, we also have $n^2-m^2<2M$, so
\begin{align*}
E_3&\ll_\eps\frac{MY^2}{T^4}\,MT^\eps\ll T^{-2\alpha+4\beta+\eps}.
\end{align*}
Since $\beta>0$, we have
$T^{-2\alpha+3\beta+\eps}\leq T^{-2\alpha+4\beta+\eps}$.
Substituting these
three estimates into the initial decomposition, together with
its $O_{\psi,\alpha,\beta}(T^{-2})$ term, proves the proposition.
\end{proof}

\section{Estimation of Long Shifts}
\label{sec:Mellin}

Throughout this section, $\omega\in\{-1,1\}$
is fixed, and $\alpha,\beta,R,M,Y$ are as in
Lemma~\ref{lem:long-local}.

In what follows, dyadic parameters run over $1,2,4,\ldots$.
Choose a fixed nonnegative function $w\in\Cc((1/2,2))$, extended
by zero to $\RR$, satisfying
\[
\sum_{N\text{ dyadic}}w(t/N)=1\qquad(t\ge1).
\]
Since $A_{m,n}=0$ for $n\le m$, we insert this partition for
each of the positive integers $m$, $n$, and $n-m$ in \eqref{eq:Hsign}.
This gives
\begin{align}
\mathcal H_\omega
={}&\frac1{TL_u}\sum_{P,Q,H\text{ dyadic}}\sum_{m<n}
\lambda_u(m)\lambda_u(n)(n/m)^{\omega\ii T}
\nonumber\\
&\qquad\times w(m/P)w(n/Q)w((n-m)/H)
A_{m,n}\cW_\psi\left(\omega\frac{n^2-m^2}{T}\right)
\nonumber\\
={}&\sum_{P,Q,H\text{ dyadic}}\mathcal U_\omega(P,Q,H),
\label{eq:H-dyadic}
\end{align}
where
\begin{equation}\label{eq:block-normalized}
\mathcal U_\omega(P,Q,H)
=\frac1{TL_u}\sum_{m,n\ge1}
\lambda_u(m)\lambda_u(n)(n/m)^{\omega\ii T}
V_{P,Q,H}^{(\omega)}(m/P,n/Q),
\end{equation}
with
\begin{align}
V_{P,Q,H}^{(\omega)}(x,y)
={}&w(x)w(y)w\left(\frac{Qy-Px}{H}\right)
\left(1-\chi\left(\frac{Qy-Px}{R}\right)\right)
\nonumber\\
&\times\chi\left(\frac{Q^2y^2-P^2x^2}{M}\right)
\cW_\psi\left(\omega\frac{Q^2y^2-P^2x^2}{T}\right).
\label{eq:Vweight}
\end{align}
Indeed, at $(x,y)=(m/P,n/Q)$ the last three factors in
\eqref{eq:Vweight} give
$A_{m,n}\cW_\psi(\omega(n^2-m^2)/T)$ for $m<n$, while the first
three give $w(m/P)w(n/Q)w((n-m)/H)$.
The factor $w((n-m)/H)$ vanishes for $n\le m$, so all index
restrictions in \eqref{eq:block-normalized} are contained in the weight.

By Lemma~\ref{lem:long-local}, a nonzero summand has
$m,n<Y$ and $R<n-m<Y$. Since $w$ is supported in $(1/2,2)$,
a triple can contribute to \eqref{eq:H-dyadic} only if
\begin{equation}\label{eq:dyadic-ranges}
1\le P,Q<2Y,\qquad R/2<H<2Y.
\end{equation}
There are
$O_{\alpha,\beta}((\log T)^3)$ such triples.

To estimate a single sum $\mathcal U_\omega(P,Q,H)$, we use
Mellin inversion. Write $V=V_{P,Q,H}^{(\omega)}$ and
\[
\widetilde V(s,z)=\int_0^\infty\int_0^\infty
V(x,y)x^{s-1}y^{z-1}\dd x\dd y.
\]

\begin{lemma}\label{lem:Mellin-norm}
For every fixed $\nu\ge0$, uniformly in $T,R,M,P,Q,H>0$
and $\omega\in\{\pm1\}$,
\begin{equation*}
\int_{\RR^2}
\left|\widetilde V\left(\tfrac12+\ii\tau_1,
\tfrac12+\ii\tau_2\right)\right|
(1+|\tau_1|+|\tau_2|)^\nu\dd\tau_1\dd\tau_2
\ll_{\psi,\nu}(1+P/H)^\nu.
\end{equation*}
\end{lemma}

\begin{proof}
Put $h=\min\{1,H/P\}$ and introduce the commuting operators
\[
D=x\partial_x+y\partial_y,\qquad E=hy\partial_y.
\]
The parameters are held fixed when differentiating.
The support of $V$, and of each $D^jE^kV$, lies in
\[
\Omega=\{(x,y)\in[1/2,2]^2:H/2\le Qy-Px\le2H\}.
\]
For fixed $y$, the allowed $x$-interval has length at most
$\min\{3/2,3H/(2P)\}$, and hence $|\Omega|\le9h/4$.

We first show that
\begin{equation}\label{eq:V-derivatives}
\|D^jE^kV\|_\infty\ll_{\psi,j,k}1
\qquad(j,k\ge0),
\end{equation}
uniformly in all the parameters. On $\Omega$,
\[
\frac{hQy}{Qy-Px}
=h\left(1+\frac{Px}{Qy-Px}\right)\le5,
\qquad
\frac{hQ^2y^2}{Q^2y^2-P^2x^2}
=\frac{hQy}{Qy-Px}\frac{Qy}{Qy+Px}\le5.
\]
For $j\ge0$ and $k\ge1$ we have
\begin{align*}
D^j(Qy-Px)&=Qy-Px,&
D^jE^k(Qy-Px)&=h^kQy,\\
D^j(Q^2y^2-P^2x^2)&=2^j(Q^2y^2-P^2x^2),&
D^jE^k(Q^2y^2-P^2x^2)&=2^{j+k}h^kQ^2y^2.
\end{align*}
Since $h\le1$, these identities give
$|D^jE^kt|\ll_{j,k}|t|$ for each of the arguments
\[
t=\frac{Qy-Px}{H},\quad
\frac{Qy-Px}{R},\quad
\frac{Q^2y^2-P^2x^2}{M},\quad
\omega\frac{Q^2y^2-P^2x^2}{T}.
\]
The chain rule therefore bounds every mixed derivative of a
factor $F(t)$ by a finite sum of $|t|^\ell|F^{(\ell)}(t)|$.
These quantities are bounded for $F=\cW_\psi$ because
$\cW_\psi$ is Schwartz, and for $F=w$ or either cutoff involving
$\chi$ because their derivatives of order at least one have fixed
compact support.
The derivatives of $w(x)$ and $w(y)$ are also uniformly bounded,
since $x,y\in[1/2,2]$ and $h\le1$. Leibniz's rule proves \eqref{eq:V-derivatives}.

Repeated integration by parts, with no boundary terms because of
compact support, gives
\[
(-s-z)^A(1+hz)^B\widetilde V(s,z)
=\int_0^\infty\int_0^\infty
D^A(1-E)^BV(x,y)x^{s-1}y^{z-1}\dd x\dd y
\]
for integers $A,B\ge0$.
On $s=1/2+\ii\tau_1$, $z=1/2+\ii\tau_2$, the integral is
$O_{\psi,A,B}(h)$ by \eqref{eq:V-derivatives} and $|\Omega|\ll h$.
Moreover,
\[
|s+z|\asymp1+|\tau_1+\tau_2|,
\qquad |1+hz|\asymp1+h|\tau_2|.
\]
Thus
\begin{equation}\label{eq:Mellin-decay}
\left|\widetilde V\left(\tfrac12+\ii\tau_1,
\tfrac12+\ii\tau_2\right)\right|
\ll_{\psi,A,B}
\frac{h}{(1+|\tau_1+\tau_2|)^A(1+h|\tau_2|)^B}.
\end{equation}

Choose integers $A,B>\nu+1$ and set
$\eta_1=\tau_1+\tau_2$, $\eta_2=h\tau_2$.
Then $\dd\tau_1\dd\tau_2=h^{-1}\dd\eta_1\dd\eta_2$ and
\[
1+|\tau_1|+|\tau_2|
\le2h^{-1}(1+|\eta_1|)(1+|\eta_2|).
\]
Integrating \eqref{eq:Mellin-decay} with the required weight yields
\begin{align*}
&\int_{\RR^2}
\left|\widetilde V\left(\tfrac12+\ii\tau_1,
\tfrac12+\ii\tau_2\right)\right|
(1+|\tau_1|+|\tau_2|)^\nu\dd\tau_1\dd\tau_2\\
&\qquad\ll_{\psi,A,B,\nu}
h^{-\nu}\int_\RR(1+|\eta_1|)^{\nu-A}\dd\eta_1
\int_\RR(1+|\eta_2|)^{\nu-B}\dd\eta_2\\
&\qquad\ll_{\psi,\nu}h^{-\nu}
\le(1+P/H)^\nu.
\end{align*}
\end{proof}

We now apply this estimate to the dyadic sums and complete the
bound for $\mathcal H_\omega$.

\begin{proposition}\label{prop:Mellin-bound}
Assume GLH.
For every $\eps>0$,
\begin{equation}\label{eq:Hbound}
\mathcal H_\omega
\ll_{\psi,\alpha,\beta,\eps}
T^{-\alpha+\beta+\eps}
\qquad(\omega=\pm1).
\end{equation}
\end{proposition}

\begin{proof}
Fix a block in \eqref{eq:H-dyadic} and abbreviate
$V=V_{P,Q,H}^{(\omega)}$. Mellin inversion in both variables
(one may see \cite[Section~3]{RS}), applied to \eqref{eq:block-normalized}, gives
\begin{equation*}
\begin{aligned}
\mathcal U_\omega(P,Q,H)
={}&\frac1{TL_u(2\pi\ii)^2}
\int_{(2)}\int_{(2)}\widetilde V(s,z)P^sQ^z\\
&\hspace{12mm}\times
L(s+\omega\ii T,u)L(z-\omega\ii T,u)\dd z\dd s.
\end{aligned}
\end{equation*}
The two Dirichlet series arise from
$(n/m)^{\omega\ii T}m^{-s}n^{-z}
=m^{-(s+\omega\ii T)}n^{-(z-\omega\ii T)}$.

For fixed parameters, $\widetilde V$ decays faster than any power
of $1+|\Im s|+|\Im z|$, uniformly for
$1/2\le\Rea s,\Rea z\le2$.
Both $L$-functions are entire and have polynomial growth in these
strips. We may therefore move one contour at a time to real part
$1/2$, without crossing a pole.

Write $s=1/2+\ii\tau_1$ and $z=1/2+\ii\tau_2$.
Since $1+T+|\tau_j\pm T|\ll T(1+|\tau_j|)$,
GLH gives, for every fixed $\nu>0$,
\[
\left|L\left(\tfrac12+\ii(\tau_1+\omega T),u\right)
L\left(\tfrac12+\ii(\tau_2-\omega T),u\right)\right|
\ll_\nu T^{2\nu}(1+|\tau_1|+|\tau_2|)^{2\nu}.
\]
As $|P^sQ^z|=\sqrt{PQ}$, Lemma~\ref{lem:Mellin-norm}
with exponent $2\nu$ yields
\begin{align*}
|\mathcal U_\omega(P,Q,H)|
&\ll_\nu\frac{\sqrt{PQ}}{TL_u}T^{2\nu}
\int_{\RR^2}
\left|\widetilde V\left(\tfrac12+\ii\tau_1,
\tfrac12+\ii\tau_2\right)\right|
(1+|\tau_1|+|\tau_2|)^{2\nu}\dd\tau_1\dd\tau_2\\
&\ll_{\psi,\nu}\frac{\sqrt{PQ}}{TL_u}
T^{2\nu}(1+P/H)^{2\nu}.
\end{align*}
For the triples in \eqref{eq:dyadic-ranges},
\[
\frac{\sqrt{PQ}}T\ll\frac YT=2T^{-\alpha+\beta},
\qquad
\frac PH\ll\frac YR=2T^{1+\beta-2\alpha}\ll T.
\]
Together with $L_u^{-1}\ll_\delta T^\delta$, this gives
\[
\mathcal U_\omega(P,Q,H)
\ll_{\psi,\alpha,\beta,\delta,\nu}
\frac{\sqrt{PQ}}T\,T^{\delta+4\nu}.
\]
For every $\eps>0$, choosing $\delta,\nu>0$ with
$\delta+4\nu<\eps$ therefore gives
\begin{equation}\label{eq:block-target}
\mathcal U_\omega(P,Q,H)
\ll_{\psi,\alpha,\beta,\eps}\frac{\sqrt{PQ}}{T}\,T^\eps,
\end{equation}
uniformly for the triples in \eqref{eq:dyadic-ranges}.
We now apply \eqref{eq:block-target} with exponent $\eps/2$ and
sum over the $O_{\alpha,\beta}((\log T)^3)$ triples in
\eqref{eq:H-dyadic}. Using $\sqrt{PQ}\ll Y$ and
$(\log(2+T))^3\ll_\eps T^{\eps/2}$, we obtain
\begin{align*}
\mathcal H_\omega
&\ll_{\psi,\alpha,\beta,\eps}
\frac{Y}{T}\,T^{\eps/2}(\log(2+T))^3\\
&\ll_{\psi,\alpha,\beta,\eps}T^{-\alpha+\beta+\eps}.
\end{align*}
This proves \eqref{eq:Hbound}.
\end{proof}

\section{Proof of Theorem~\ref{thm:main} and  Corollary~\ref{cor:interval}}

\begin{proof}[Proof of Theorem~\ref{thm:main}]
Combining \eqref{eq:mass-decomposition} and
\eqref{eq:short-decomposition} with
Lemmas~\ref{lem:diagonal-bound}, \ref{lem:short-period-bound}, and
\ref{lem:short-correction-bound}, and
Propositions~\ref{prop:long-reduction} and~\ref{prop:Mellin-bound}
gives, for every fixed $\eta,\sigma>0$,
\begin{align}
I_u(\psi)-M(\psi)
\ll_{\psi,\alpha,\beta,\eta,\sigma}\;&
T^{-1/2+2\eta+\alpha}
+T^{-1+2\alpha+\sigma}
+T^{-\alpha+\beta+\sigma}
\nonumber\\
&+T^{-1+\beta+\sigma}
+T^{-2\alpha+4\beta+\sigma}
+T^{-2}.
\label{eq:all-errors}
\end{align}
Take $\alpha=1/4$ and $\beta=\eta=\sigma=\eps$, where
$0<\eps<1/16$. These choices satisfy the parameter conditions in Lemma~\ref{lem:long-local},
and \eqref{eq:all-errors} gives the following estimate, which proves
\eqref{eq:main} after relabeling $2\eps$ as $\eps$:
\[
I_u(\psi)-M(\psi)\ll_{\psi,\eps}T^{-1/4+2\eps}.\qedhere
\]
\end{proof}

\begin{proof}[Proof of Corollary~\ref{cor:interval}]
Write $A=\log a$, $B=\log b$, and fix $\delta>0$.
Choose $d>0$ with $d<(B-A)/4$ and $4d<\delta$.
Let $\eta_d\in\Cc((-d,d))$ be nonnegative with integral one, and set
\[
\phi_-=\mathbf1_{[A+d,B-d]}*\eta_d,\qquad
\phi_+=\mathbf1_{[A-d,B+d]}*\eta_d,\qquad
\psi_\pm(y)=\phi_\pm(\log y).
\]
These are smooth compactly supported functions on $(0,\infty)$.
Convolution with the nonnegative kernel gives
\[
0\le\psi_-\le\mathbf1_{[a,b]}\le\psi_+,
\qquad
\int_0^\infty(\psi_+-\psi_-)\frac{\dd y}{y}=4d<\delta.
\]
The pointwise inequalities also hold at the endpoints, which
have zero measure in all the integrals.
Positivity of $|u(iy)|^2$ bounds the interval integral
between $I_u(\psi_-)$ and $I_u(\psi_+)$.
For these fixed functions, let $T$ tend to infinity using
Theorem~\ref{thm:main} with $\eps=1/12$. Their limiting integrals approach
$(6/\pi)(B-A)$ as $d\downarrow0$.
This proves \eqref{eq:interval}.
\end{proof}

\end{document}